\documentclass[11pt]{article}

\usepackage{tikz}
\usepackage{amsmath, amssymb, amsthm}
\usepackage{t1enc}
\usepackage[utf8]{inputenc}
\usepackage[english]{babel}
\usepackage[a4paper]{geometry}

\title{General bounds for metric dimension of diameter-two graphs with applications}

\date{}

\author{\'Akos Beke\thanks{Department of Mathematics, ELTE E\"otv\"os Lor\'and University,
		Savaria University Centre, Károlyi Gáspár tér 4., 9700 Szombathely, Hungary,
		e-mail address: beke.akos@sek.elte.hu}}

\newtheorem{theorem}{Theorem}[section]
\newtheorem{lemma}[theorem]{Lemma}

\newtheorem{proposition}[theorem]{Proposition}

\theoremstyle{definition}
\newtheorem{definition}[theorem]{Definition}
\newtheorem{notation}[theorem]{Notation}

\newtheorem{remark}[theorem]{Remark}

\newcommand{\N}[0]{\mathbb{N}}

\newcommand{\R}[0]{\mathbb{R}}

\newcommand{\E}[0]{\mathbb{E}}
\newcommand{\F}[0]{\mathbb{F}}

\newcommand{\PG}[0]{\mathrm{PG}}
\newcommand{\SRG}[0]{\mathrm{SRG}}

\begin{document}
	
\maketitle

\begin{abstract}
	We study the metric dimension of several classes of graphs with diameter two. First,
	we establish a general lower bound for the metric dimension of graphs with diameter
	two in terms of their order and maximum degree. We then apply our results to
	Grassmann graphs, Moore graphs, and bilinear forms graphs. For Kneser graphs,
	we obtain new upper bounds, while for $G_q(n,2)$ Grassmann graphs we establish new upper 
	bounds provided by constructions of resolving sets. For Moore graphs of diameter
	two, we derive improved lower bounds, and for $H_q(n,2)$ bilinear forms graphs, we obtain new
	upper and lower bounds on their metric dimension.
\end{abstract}

\section{Preliminaries}

\begin{definition} Let $G=(V, E)$ be a graph, $S\subseteq V$. Let us fix an order
	of the elements of $S=(s_1, s_2, \dots, s_m)$. The
	\emph{distance vector of $v\in V$ with respect to $S$} is
	\[R(v|S)=(d(v,s_1), d(v,s_2), \dots, d(v,s_m)).\]
	We say that a set $W\subseteq V$ is \emph{resolved by the set $S$} 
	if $R(v|S)\ne R(u|S)$ for any $v\ne u\in W$, that is, there is a vertex
	$s\in S$ such that $d(v,s)\ne d(u,s)$.

	$S$ is called a \emph{resolving set} of $G$ if it resolves the set $V$.
	A smallest resolving set is called \emph{optimal}, its size is called the
	\emph{metric dimension of the graph} and it is denoted by $\mu(G)$.
\end{definition}
	
The concept of metric dimension can be discussed in any metric
space, and it already appeared in 1953  \cite{AppOfDisGeo}.
The concept of resolving sets was introduced by Slater in 1975 \cite{ResSetSlater} and
independently by Harary and Melter in 1976 \cite{ResSetHararyMelter}.
Many results of the topic have been gathered in \cite{BasMetDim}, \cite{ResSetMetDim} and
\cite{SurveyMetDim2023}. Since then, the metric dimension of various graph classes has been
studied, including numerous graphs arising from finite geometries
\cite{Bailey,KnesJohnMetDim,GrassmannUpperPaper,MetDimAffBiaffGenQuad,ResSetHigDimProjPL,ResSetHigDimProj,LocNumMetDimGrDiam2,HegTakResSet,HegSzilTakMetSmallOrd}.

In this paper we investigate resolving sets of Kneser, Grassmann, Moore and Bilinear forms graphs,
obtaining new results for both classes. These graph families have previously been studied
\cite{KnesJohnMetDim,GrassmannUpperPaper,LocNumMetDimGrDiam2,FengWangBilFormGrMetric}.

\begin{notation}
	Let $f,g:\R^n\to\R$ be two multivariable functions. We say that $f(\vec x)$
	\emph{asymptotically equals} $g(\vec x)$ and we denote it by $f(\vec x)\sim g(\vec x)$, if
	\[\lim_{\vec x\to\infty}\frac{f(\vec x)}{g(\vec x)}=1,\]
	where $\vec x\to\infty$ means that all variables tend to infinity independently, along any
	possible path in $\R^n$. In other words, the quotient $\frac {f(\vec x)}{g(\vec x)}$ must
	converge to $1$ regardless of how the variables tend to infinity.
\end{notation}

\section{General results}

\subsection{Lower bounds}

We will use the following known combinatorial inequality:
\begin{proposition}[\cite{StaJukExtBook} Proposition 1.4.]\label{InequalityBinomialSumProp}
	If $u\in\N$, and $0<u\le k$, then
	\begin{equation}\label{InequalityBinomialSum}
		\sum_{i=0}^u\binom{k}{i}<\left(\frac{ek}{u}\right)^u.
	\end{equation}
\end{proposition}

\begin{lemma}\label{LnuPerLnVMinusLnKPlusOne}
	Let $k, v$ be positive integers with $k\le v$. If $v\ge10$ or $k\le\frac{ev}4$ and $v\ge4$, then 
	\[\frac{\ln v}{\ln v-\ln k +1}\le\log_2(v)-1.\]
\end{lemma}

\begin{proof}
	By reformulating
	\begin{equation}\label{uwhb1fjfniyqgunrvjhuwIJUDVUBR}
	\frac{\ln v}{\ln\frac{ev}k}\le\frac{\ln v-\ln 2}{\ln2}.
	\end{equation}
	If $k\le\frac{ev}4$ this follows from
	\[\frac{\ln v}{\ln4}\le\frac{\ln v-\ln 2}{\ln2}.\]
	By rearranging
	\[\ln2\ln v\le\ln4\ln v-\ln2\ln4=2\ln2\ln v-\ln2\ln4,\]
	or equivalently
	\[\ln2\ln4\le\ln2\ln v,\]
	and this is equivalent to the assumption $v\ge4$. Thus the inequality holds if $v\ge4$ and
	$k\le\frac{ev}4$.
	
	Since $\frac vk\ge1$, Formula (\ref{uwhb1fjfniyqgunrvjhuwIJUDVUBR}) follows from
	\[\ln v\le\frac{\ln v-\ln 2}{\ln2},\]
	By rearranging
	\[\frac{\ln2}{1-\ln2}\le\ln v.\]
	The left hand side is approximately $2.259$ and the right hand side is increasing;
	$\ln9\approx2.197$ and $\ln10\approx2.3$, therefore the inequality holds if $v\ge10$.
\end{proof}

There is a straightforward general lower bound for any graph $G$ with diameter $d$ on $n$
vertices: 
\begin{equation}\label{ChartrandsLowerBoundInequalityForm}
	d^{\mu(G)}+\mu(G)\geq n
\end{equation}
(see e.g.\ \cite{ResSetMetDim} Theorem 1.).

\begin{theorem} Let  $G=(V,E)$ be a graph with diameter 2 and $v=|V|$. Then
	\begin{equation}\label{ChartrandsLowerBound}
		\mu(G)>\log_2(v)-1=\frac{\ln(v)}{\ln2}-1.	
	\end{equation}
\end{theorem}
\begin{proof}
	From inequality (\ref{ChartrandsLowerBoundInequalityForm}), we get
	\[2^{\mu(G)+1} > 2^{\mu(G)}+\mu(G)\ge v,\]
	Taking logarithms proves the strict inequality.
\end{proof}

The following proposition modifies the author's earlier result, Proposition 16 in
\cite{ResSetMetDimSymDes}, to the setting of non-bipartite graphs of diameter two. The
proof employs essentially the same method as in \cite{ResSetMetDimSymDes}, where the result
was established for bipartite graphs of diameter three.

\begin{lemma}\label{MainPropositionGeneralGraphs}
	Let $G=(V,E)$ be a graph with maximum degree $k$ and diameter 2,
	$v=|V|$, let $u$ be a positive integer and let $s\in(0,1)$. If $S$ is a resolving set and
	$u\le |S|$, then
	\begin{equation}\label{MainPropositionGeneralGraphsFormula}
	|S|\ge
	\min\left\{\frac{sv(u+1)}{k}\ ,
	\ \frac{u\sqrt[u]{(1-s)v}}{e\sqrt[u]{2}}\ ,
	\ \frac{(1-s)v}2
	\right\}.
\end{equation}
\end{lemma}

\begin{proof}
	\textbf{Case 1:} Suppose that there are at least $sv$ vertices in $V\setminus S$ such
	that each of them has at least $u+1$ neighbors in $S$.
	By double counting the set
	\[\{(x,y)\in S\times (V\setminus S)\ :\ \{x,y\}\in E\},\]
	we get
	\[|S|k\ge sv(u+1).\]

	\textbf{Case 2:} There are less than $sv$ vertices in $V\setminus S$ which have at least $u+1$
	neighbors in $S$. Therefore, there are at least $(1-s)v-|S|$ vertices in $V\setminus S$
	which have at most $u$ neighbors in $S$. Let
	\[V_2=\{w\in V\setminus S\ :\ |N(w)\cap S|\le u\}.\]
	If $w_1\ne w_2\in V_2$, then $R(w_1|S)\ne R(w_2|S)$.
	The number of different sequences is at least $(1-s)v-|S|$ and the number of
	sequences is at most $\sum_{i=0}^u\binom{|S|}{i}$, as there can be at most $u$ ones in
	each sequence. Therefore
	\begin{equation}\label{Case2SumFigureGeneral}
		(1-s)v-|S|\le\sum_{i=0}^u\binom{|S|}{i}\le\left(\frac{e|S|}{u}\right)^u.
	\end{equation}
	By rearrangement
	\begin{equation}\label{FigureGeneralBeforeSubcases}
	|S|^u+\frac{u^u}{e^u}|S|\ge\frac{u^u}{e^u}(1-s)v.
	\end{equation}
	\textbf{Case 2.1:}  $|S|^u\ge\frac{u^u}{e^u}|S|$. Then $2|S|^u\ge\frac{u^u}{e^u}(1-s)v$, and,
	by rearranging
	\[|S|\ge\frac{u\sqrt[u]{(1-s)v}}{e\sqrt[u]{2}}.\]
	\textbf{Case 2.2:}  $|S|^u\le\frac{u^u}{e^u}|S|$. Then
	$2\frac{u^u}{e^u}|S|\ge\frac{u^u}{e^u}(1-s)v$. Thus
	\[|S|\ge\frac{(1-s)v}2.\]

	As one of the three cases holds, we have proved the assertions.
\end{proof}

\begin{theorem}\label{SimpleLowerBoundDiamTwo}
Let $G=(V,E)$ be a graph of maximum degree $k$ with diameter 2,
$v=|V|$. Then
\[\mu(G)\ge\frac{2(v-1)}{k+4}.\]
\end{theorem}

\begin{proof}
In the proof of Lemma \ref{MainPropositionGeneralGraphs} let $u=1$. Then, in the second
case, Formula (\ref{Case2SumFigureGeneral}) can be modified by not using Proposition
\ref{InequalityBinomialSumProp}:
\[(1-s)v-|S|\le\sum_{i=0}^1\binom{|S|}{i}=1+|S|.\]
By rearrangement
\[|S|\ge\frac{(1-s)v-1}2.\]

With this modification Lemma \ref{MainPropositionGeneralGraphs} claims that
\[|S|\ge
\min\left\{\frac{2sv}{k}\ ,
\ \frac{(1-s)v-1}2
\right\}.
\]
The two formulas are equal to each other if $s=\frac{k(v-1)}{v(k+4)}$. By substituting
\[|S|\ge\frac{2(v-1)}{k+4}.\]
\end{proof}

\begin{theorem}\label{MainTheoremMetricDimensionDiameterTwo}
Let $G=(V,E)$ be a graph with maximum degree $k\ge3$ and diameter 2,
$v=|V|$. If $v\ge10$ or $k\le\frac{ev}4$ and $v\ge4$, then
\[\mu(G)\ge\frac{uv}{k\sqrt[u]{2u}},\]
where $u=\left\lfloor\frac{\ln v}{\ln v-\ln k+1} \right\rfloor$.
\end{theorem}

\begin{proof}
By Lemma \ref{MainPropositionGeneralGraphs}
\[\mu(G)\ge\min\left\{\frac{sv(u+1)}{k}\ ,
\ \frac{u\sqrt[u]{(1-s)v}}{e\sqrt[u]{2}}\ ,
\ \frac{(1-s)v}2
\right\}.\]
Let $u=\left\lfloor \frac{\ln v}{\ln v-\ln k+1} \right\rfloor$. The inequality
$\frac{\ln v}{\ln v-\ln k+1}\ge1$ is equivalent to $k\ge e$, which is implied by $k\ge3$.
And if $u=1$, then by the preceding theorem
\[\mu(G)\ge\frac{2(v-1)}{k+4}\ge\frac{uv}{k\sqrt[u]{2u}}=\frac v{2k}.\]
The inequality $\frac{2(v-1)}{k+4}\ge\frac v{2k}$ is equivalent to $k\ge\frac{4v}{3v-4}$,
which holds, since $v\ge4$.

From now on, we assume that $u\ge2$.
By (\ref{ChartrandsLowerBound}) and the Lemma \ref{LnuPerLnVMinusLnKPlusOne} 
$u=\left\lfloor \frac{\ln v}{\ln v-\ln k+1} \right\rfloor$ can be used as parameter $u$ in
Lemma \ref{MainPropositionGeneralGraphs}.
\[u\le\frac{\ln v}{\ln v-\ln k+1}=\frac{\ln v}{\ln\frac{ev}k}=\log_{\frac{ev}k}v,\]
thus $\left(\frac{ev}k\right)^u\le v$ or equivalently $\frac{\sqrt[u]{v}}e\ge\frac{v}k$,
and therefore
\[\mu(G)\ge\min\left\{\frac{sv(u+1)}{k}\ ,
\ \frac{u\sqrt[u]v}e \sqrt[u]{\frac{1-s}{2} }\ ,
\ \frac{(1-s)\left(\frac{ev}k\right)^u}2
\right\}\]
\[\ge\min\left\{\frac{suv}{k}\ ,
\ \frac{uv}k \sqrt[u]{\frac{1-s}{2} }\ ,
\ \frac {uv}k\frac{(1-s)\left(\frac vk\right)^{u-1}e^u}{2u}
\right\}\ge\]\[
\frac{uv}k\min\left\{s\ ,
\  \sqrt[u]{\frac{1-s}{2} }\ ,
\ \frac{(1-s)e^u}{2u}
\right\}.\]
Let $s=\frac{u-1}{u}$. Then
\[\mu(G)\ge\frac{uv}k\min\left\{\frac{u-1}{u}\ ,
\  \sqrt[u]{\frac{1}{2u} }\ ,
\ \frac{e^u}{2u^2}
\right\}.\]
It is easy to see that $\frac{u-1}{u}\le\frac{e^u}{2u^2}$. So we only have to prove that
\[\sqrt[u]{\frac{1}{2u}}\le\frac{u-1}{u}.\]
This holds for $u=2$. So we can assume that $u\ge3$. By taking the logarithm of both sides
\[\frac1u\ln\frac1{2u}\le\ln\left(1-\frac1u\right).\]
By the known inequality $\ln(1-x)\ge\frac{x}{x-1}$ ($x<1$), this follows from
\[\frac1u\ln\frac1{2u}\le\frac1{1-u}.\]
Or equivalently
\begin{equation}\label{kjdbjhsbfrhkahdfljhkbgflklkfsdhb}
\ln(2u)\ge\frac u{u-1}.
\end{equation}
Let us calculate the derivatives of both sides:
\[\frac d{du}\ln(2u)=\frac1u>0.\]
Thus the left side of (\ref{kjdbjhsbfrhkahdfljhkbgflklkfsdhb})
is increasing in $u$, and the right side is decreasing:
\[\frac d{du}\frac u{u-1}=\frac{-1}{(u-1)^2}<0.\]
Since inequality (\ref{kjdbjhsbfrhkahdfljhkbgflklkfsdhb}) holds for $u=3$, it also holds for $u\ge3$.
\end{proof}

\subsection{Upper bounds}


\begin{definition}
	Let $G=(V,E)$ be a graph that is neither complete nor empty. $G$ is called
	\emph{strongly regular} with \emph{parameters $(v, k, a, b)$} if 
	\begin{itemize}
		\item $v=|V|$,
		\item $G$ is $k$-regular,
		\item every two adjacent vertices have $a$ common neighbors, and
		\item every two non-adjacent vertices have $b$ common neighbors.
	\end{itemize}
	We denote a strongly regular graph with parameters $(v, k, a, b)$ by $\SRG(v,k,a,b)$.
	A strongly regular graph is called \emph{primitive} if both itself and its complement
	are connected, otherwise \emph{imprimitive}.
\end{definition}

It is easy to see that the diameter of a primitive strongly regular graph is two, thus our
theorems can be used to obtain a lower
bound to the metric dimension. As for the upper bounds, the following theorem
is a corollary of Theorem 3.22 in \cite{BasMetDim} which holds for primitive
distance regular graphs, and so it also holds for primitive strongly regular graphs. It was
originally given by Babai in \cite{BabaiPermutation}.

\begin{theorem}\label{BabaiUpperBoundStrongly}
	Let $G=(V,E)$ be a primitive strongly regular graph, $v=|V|$, and let $k$ be the degree of a vertex.
	Then
	\[\mu(G)<4\frac v{v-M}\ln v,\]
	where $M=\max\{k, v-k-1\}$.
\end{theorem}

We establish a new upper bound on the metric dimension of regular graphs, which also
applies to strongly regular graphs.

\begin{theorem}\label{UpperBoundAlmostStronglyRegular}
	Let $G=(V,E)$ be a $k$-regular graph with $v=|V|$ and $k\ge1$.
	Let $a$ and $b$ be the smallest non-negative integers such that every two adjacent
	vertices have at most $a$ common neighbors, and every two non-adjacent
	vertices have at most $b$ common neighbors (for example $\SRG(v,k,a,b)$),
	and let $N=\max\{a,b-1\}$, and $v>2k-2N$. Then
	\[\mu(G)\le\left\lceil\frac{v\ln\frac v{\sqrt 2}}{k-N}\right\rceil.\]
\end{theorem}

\begin{proof}
	First we compute the probability that a randomly chosen vertex $s$ does not resolve
	$\{x,y\}$.
	The set $\{x,y\}$ is resolved by $x$, $y$ or any neighbor of $x$ or $y$ that is not a
	common neighbor.
	If $x$ and $y$ are two adjacent vertices, then
	\[\Pr(s \textrm{ does not resolve } \{x,y\})=1-\Pr(s \textrm{ resolves } \{x,y\})\]
	\[\le1-\frac{2k-2a}v=\frac{v-2k+2a}v.\]
	If $x$ and $y$ are two non-adjacent vertices, then
	\[\Pr(s \textrm{ does not resolve } \{x,y\})=1-\Pr(s \textrm{ resolves } \{x,y\})\]
	\[\le1-\frac{2k-2b+2}v=\frac{v-2k+2b-2}v.\]
	Choose $s_1, s_2, \dots, s_m$ independently with replacement, and let
	$S=\{s_1, s_2, \dots, s_m\}$ be the set of chosen vertices.
	Let
	\[I_{\{x,y\}}=
	\begin{cases}
		1, & \text{if $S$ does not resolve  } \{x,y\},\\
		0, & \text{otherwise}.
	\end{cases}\]

	\[\E(I_{\{x,y\}})\le
	\begin{cases}
		\left(\frac{v-2k+2a}v\right)^m, & \text{if } x,y \text{ are adjacent},\\
		\left(\frac{v-2k+2b-2}v\right)^m, & \text{otherwise}.
	\end{cases}\]
	
	Let $X$ denote the number of pairs $\{x,y\}$ that are not resolved by $S$.
	\[\E(X)\le \sum_{\substack{\text{adjacent} \\ \text{pairs}}}\left(\frac{v-2k+2a}v\right)^m+
	\sum_{\substack{\text{non-adjacent} \\ \text{pairs}}}\left(\frac{v-2k+2b-2}v\right)^m\]
	\[\le\frac{vk}2 \left(\frac{v-2k+2a}v\right)^m+
	\frac{v(v-k)}2\left(\frac{v-2k+2b-2}v\right)^m \]
	\[\le\frac{vk}2 \left(\frac{v-2k+2N}v\right)^m+
	\frac{v(v-k)}2\left(\frac{v-2k+2N}v\right)^m=\frac{v^2}2 \left(\frac{v-2k+2N}v\right)^m.\]
	If $\E(X)<1$ then there exists a resolving set with cardinality $m$ and thus $\mu(G)\le m$.
	The inequality
	\[\frac{v^2}2 \left(\frac{v-2k+2N}v\right)^m<1\]
	is equivalent to
	\[m>\frac{2\ln\frac v{\sqrt 2}}{\ln\frac v{v-2k+2N}},\]
	Using the known bound $\ln x>1-\frac1x$, where $x>1$
	\[m\ge\frac{2\ln\frac v{\sqrt 2}}{1-\frac{v-2k+2N}v}=
	\frac{2v\ln\frac v{\sqrt 2}}{2k-2N},\]
	therefore
	\[\mu(G)\le\left\lceil\frac{v\ln\frac v{\sqrt 2}}{k-N}\right\rceil.\]
\end{proof}

In \cite{BabaiComplexity}, Babai proved a similar upper bound:
\begin{theorem}[\cite{BabaiComplexity}, Lemma 3.2]
Using the notation of Theorem \ref{UpperBoundAlmostStronglyRegular},
\[\mu(G)\le\left\lceil\frac{v\ln v}{k-N}\right\rceil.\]
\end{theorem}

Comparing the two bounds, we see that our result improves upon Babai's bound by an additive term of
approximately $\frac{v\ln\sqrt2}{k-N}$.

\section{Kneser graphs}

\begin{definition}
	The $K(n,k)$ \emph{Kneser graph} is a graph on vertex set $V$ that is the set of $k$-element
	subsets of $[n]=\{1, 2, \dots, n\}$ and two vertices are adjacent if they are disjoint
	subsets of $[n]$.
\end{definition}

Bailey et al. proved the following upper bounds.
\begin{theorem}[\cite{KnesJohnMetDim} Theorem 7. and Theorem 8.]
	If $n>2k$, then
	\begin{equation}\label{BaileyUpperboundNGE2}
	\mu(K(n,k))\le\left\lceil\frac n{2k-1}\right\rceil\left({{2k-1}\choose k}-1\right)
	\sim\frac{n4^{k-1}}{k\sqrt{k\pi}}.
	\end{equation}
	If $\lfloor5k/2 \rfloor\le n\le3k-2$, then
	\[\mu(K(n,k))\le2\binom{n-k}k.\]
\end{theorem}

Bonato et al. proved the following bounds.

\begin{theorem}[\cite{LocNumMetDimGrDiam2} Corollary 2.4. and Lemma 2.8.]
	If $n\ge3k$ and $k\ge3$, then
	\[\mu(K(n,k))\ge\frac n2+\frac nk.\]

	For every integer $k\ge3$, let $H(k,5)$ be a $k$-uniform $\lceil k/2+1\rceil$-regular
	hypergraph of girth 5, and let $m$ be the number of vertices in $H(k,5)$. Then
	$m\le w2^{w(w-1)/2}$, where $w=k^{\lceil k/2+1\rceil}$. Let $n\ge3k$.

	For even integers $k\ge4$,
	\[\mu(K(n,k))
	\le\frac n2+\frac nk+\left(\frac m2+\frac mk\right)\left\lceil\frac{n'}m\right\rceil,\]
	and for odd integers $k\ge3$,
	\[\mu(K(n,k))
	\le\frac n2+\frac nk+\frac n{2k}+\left(\frac m2+\frac mk+\frac m{2k}\right)
	\left\lceil\frac{n'}m\right\rceil,\]
	where $n'$ is the smallest non-negative integer such that $n'\equiv n\pmod m$.
\end{theorem}

\begin{remark}
	It is worth noting that, when $\left\lceil\frac{n'}m\right\rceil\ne0$, that is $m\nmid n$,
	the upper bounds in the above theorem contain
	an additional term proportional to $m$. Since $m$ is bounded only by $w2^{w(w-1)/2}$,
	where $w=k^{\lceil k/2+1\rceil}$, this term can be extremely large, even for relatively
	small values of $k$. Consequently, although the above bounds are asymptotically strong,
	their practical usefulness may be limited when $m\nmid n$.
\end{remark}

\begin{theorem}\label{UpperBoundKneserGraphMetricDimension}
	If $k\ge2$ and $n\ge2k+1$, then
	\[\mu(K(n,k))\le\left\lceil\frac{{n \choose k}\ln\frac {{n \choose k}}{\sqrt2}}
	{{{n-k-1} \choose {k-1}}}\right\rceil 
	\le \left\lceil n\left(\frac{n-k+1}{n-2k+1}\right)^{k-1}\ln\frac{en}k\right\rceil.\]
\end{theorem}
\begin{proof}
If $x,y\subseteq[n]$ are two adjacent vertices, then they have ${{n-2k} \choose k}$ common neighbors.
If $x,y\subseteq[n]$ are two non-adjacent vertices, $z=|x\cap y|$, then they have
${{n-2k+z} \choose k}$ common neighbors. The larger $z$ is, the more common neighbors they have,
and $z$ is at most $k-1$. Thus two non-adjacent vertices have at most ${{n-k-1} \choose k}$
common neighbors.

Let $N={{n-k-1} \choose k}>\max\left\{{{n-2k} \choose k},{{n-k-1} \choose k}-1\right\}.$
By Theorem \ref{UpperBoundAlmostStronglyRegular}
\[\mu(K(n,k))\le \left\lceil\frac{v\ln\frac v{\sqrt2}}{d-N}\right\rceil
=\left\lceil\frac{{n \choose k}\ln\frac {{n \choose k}}{\sqrt2}}{{{n-k} \choose k}-{{n-k-1} \choose k}}\right\rceil
=\left\lceil\frac{{n \choose k}\ln\frac {{n \choose k}}{\sqrt2}}{{{n-k-1} \choose {k-1}}}\right\rceil.\]
Now we give an upper bound for $\frac{{n \choose k}}{{{n-k-1} \choose {k-1}}}$:
\[\frac{{n \choose k}}{{{n-k-1} \choose {k-1}}}=\frac{n!(k-1)!(n-2k)!}{k!(n-k)!(n-k-1)!}
=\frac nk\frac{(n-1)(n-2)\dots(n-(k-1))}{(n-k-1)(n-k-2)\dots(n-k-(k-1))}\]
\[=\frac nk\prod_{i=1}^{k-1}\frac{n-i}{n-k-i}\le\frac nk\left(\frac{n-k+1}{n-2k+1}\right)^{k-1}.\]
We also use the known upper bound ${n \choose k}<\left(\frac{en}{k}\right)^k$, and therefore
\[\mu(K(n,k))\le\left\lceil n\left(\frac{n-k+1}{n-2k+1}\right)^{k-1}\ln\frac{en}k\right\rceil.\]
The assumption $v>2d-2N$ also holds: it is equivalent to $\frac vd>2-\frac{2N}d$. By substitution,
we have to prove that
\[\frac{{n \choose k}}{{{n-k} \choose k}}>2-2\frac{{{n-k-1} \choose k}}{{{n-k} \choose k}}.\]
Since
\[\frac{{n \choose k}}{{{n-k} \choose k}}=\prod_{i=0}^{k-1}\frac{n-i}{n-k-i}\ge\frac n{n-k},\]
and
\[\frac{{{n-k-1} \choose k}}{{{n-k} \choose k}}=\frac{(n-k-1)!k!(n-2k)!}{k!(n-2k-1)!(n-k)!}=\frac{n-2k}{n-k}\]
the desired inequality follows from
\[\frac n{n-k}>2-2\frac{n-2k}{n-k},\]
which is equivalent to $n>2k$.
\end{proof}

\begin{remark}
	Let $c=\frac nk$. Then $\frac{n-k+1}{n-2k+1}=1+\frac k{ck-2k+1}\le1+\frac1{c-2}$. Thus, our upper
	bound is at most $n\left(1+\frac1{c-2}\right)^{k-1}\ln\frac{en}k$. For $c\ge3$, that is, $n\ge3k$,
	we have $\frac{n-k+1}{n-2k+1}\le2$, and hence our upper bound is at most $n2^{k-1}\ln\frac{en}k$.
	Bailey's upper bound in Formula (\ref{BaileyUpperboundNGE2}) is asymptotically
	$\frac{n4^{k-1}}{k\sqrt{k\pi}}$. Therefore, a rough comparison suggests that our bound can be
	substantially better when roughly $n\lesssim\frac ke\exp\frac{2^{k-1}}{k\sqrt{k\pi}}$,
	although Bailey's bound is asymptotically better if $n\gtrsim\frac ke\exp\frac{2^{k-1}}{k\sqrt{k\pi}}$.
\end{remark}

Bailey et al. proved the following upper bounds, too:
\begin{theorem}[\cite{KnesJohnMetDim} Corollary 14. and Corollary 15.]
	If $q$ is a prime power then
	\begin{enumerate}
		\item if $q\ge3$, then $\mu(K(q^2,q))\le q^2+q$,
		\item $\mu(K(q^3,q))\le q^2(q+2)$,
		\item $\mu(K((q+1)(q^3+1),q+1))\le (q^2+1)(q^3+1)$,
		\item $\mu(K((q^2+1)(q^5+1),q^2+1))\le (q^3+1)(q^5+1)$.
\end{enumerate}
\end{theorem}
Theorem \ref{UpperBoundKneserGraphMetricDimension} provides better
upper bounds asymptotically for the third and fourth cases:
\begin{theorem}
	Let $m\ge2$ be an integer. Then
\[\mu(K((m+1)(m^3+1),m+1))\]
\[\le \left\lceil(m+1)(m^3+1)\left(\frac{m^4+m^3+1}{m^4+m^3-m}\right)^m
\ln (e(m^3+1))\right\rceil\sim 3m^4\ln m,\]
\[\mu(K((m^2+1)(m^5+1),m^2+1))\le\]
\[\left\lceil(m^2+1)(m^5+1)\left(\frac{m^7+m^5+1}{m^7+m^5-m^2}\right)^{m^2}\ln(e(m^5+1))\right\rceil
\sim 5m^7\ln m.\]
\end{theorem}

\section{Grassmann graphs}
\begin{definition}
	The $G_q(n,k)$ \emph{Grassmann graph} is a graph on vertex set $V$ that is the set of $k$-dimensional
	subspaces of the vector space $\F_q^n$ and two vertices are adjacent if they intersect in a
	$(k-1)$-dimensional subspace.
\end{definition}

In this section we consider the $G_q(n,2)$ Grassmann graphs with $n\ge4$, that is, their vertices are the lines
of the projective space $\PG(n-1,q)$, and two vertices are adjacent if they have one common point.
The $G_q(n,2)$ Grassmann graph is strongly regular with parameters
\[\left(\begin{bmatrix} n \\ 2 \end{bmatrix}_q\ ,
\ (q+1)\left(\begin{bmatrix} n-1 \\ 1 \end{bmatrix}_q-1\right)\ ,
\ \begin{bmatrix} n-1 \\ 1 \end{bmatrix}_q+q^2-2\ ,\ (q+1)^2\right).\]

Let $2\le k\le\frac n2$.
In \cite{GrassmannUpperPaper} Bailey and Meagher proved that 
\begin{equation}\label{BaiMeaGrasUpp}
\mu(G_q(n,k))\le\begin{bmatrix} n \\ 1 \end{bmatrix}_q=\frac{q^n-1}{q-1}.
\end{equation}

In \cite{GrassmannUpperPaper} Bailey and Meagher mentioned two other upper bounds given
by Babai in \cite{BabaiPermutation}:
\begin{equation}\label{BabaiGrasFirstUpp}
\mu(G_q(n,k))<4\sqrt{\begin{bmatrix} n \\ k \end{bmatrix}_q}
\ln\begin{bmatrix} n \\ k \end{bmatrix}_q,
\end{equation}
and by Theorem 3.22 of \cite{BasMetDim}
\begin{equation}\label{BabaiGrasSecondUpp}
\mu(G_q(n,k))<2k\frac{\begin{bmatrix} n \\ k \end{bmatrix}_q}{\begin{bmatrix} n \\ k \end{bmatrix}_q-M}
\ln\begin{bmatrix} n \\ k \end{bmatrix}_q,
\end{equation}
where
\[M=\max_{0\le j\le k}q^{j^2}\begin{bmatrix} n-k \\ j \end{bmatrix}_q
\begin{bmatrix} k \\ j \end{bmatrix}_q.\]
Whereas \eqref{BaiMeaGrasUpp} does not depend on $k$, both of Babai's upper bounds do.
Let us evaluate them at $k=2$. Formula \eqref{BabaiGrasFirstUpp} is
\[\mu(G_q(n,2))<4\sqrt{\begin{bmatrix} n \\ 2 \end{bmatrix}_q}
\ln\begin{bmatrix} n \\ 2 \end{bmatrix}_q=\]
\[4\sqrt{\frac{(q^n-1)(q^{n-1}-1)}{(q^2-1)(q-1)}}
\ln\frac{(q^n-1)(q^{n-1}-1)}{(q^2-1)(q-1)}\sim4q^{n-2}(2n-4)\ln q.\]
To evaluate \eqref{BabaiGrasSecondUpp}, first we have to calculate $M$:
If $j=0$, then
\[q^{j^2}\begin{bmatrix} n-k \\ j \end{bmatrix}_q
\begin{bmatrix} k \\ j \end{bmatrix}_q=1.\]
If $j=1$, then
\begin{equation}\label{LargestMForSmallQAndSmallNBabai}
q^{j^2}\begin{bmatrix} n-k \\ j \end{bmatrix}_q
\begin{bmatrix} k \\ j \end{bmatrix}_q=q\frac{q^{n-2}-1}{q-1}\ \frac{q^2-1}{q-1}=
q(q+1)\frac{q^{n-2}-1}{q-1}.
\end{equation}
If $j=2$, then
\begin{equation}\label{LargestMBabai}
q^{j^2}\begin{bmatrix} n-k \\ j \end{bmatrix}_q
\begin{bmatrix} k \\ j \end{bmatrix}_q=
q^4\frac{(q^{n-2}-1)(q^{n-3}-1)}{(q^2-1)(q-1)}.
\end{equation}
It is easy to see that the last one is the greatest value of the formula if $q\ge3$ and
$n\ge4$ or $q\ge2$ and $n\ge5$, thus $M$ is equal to Formula \eqref{LargestMBabai}, and therefore,
for $k=2$, Formula \eqref{BabaiGrasSecondUpp} yields
\[\mu(G_q(n,2))<4\frac{\begin{bmatrix} n \\ 2 \end{bmatrix}_q}{\begin{bmatrix} n \\ 2 \end{bmatrix}_q-M}
\ln\begin{bmatrix} n \\ 2 \end{bmatrix}_q
=4\frac{\frac{(q^n-1)(q^{n-1}-1)}{(q^2-1)(q-1)}}
{\frac{(q^n-1)(q^{n-1}-1)}{(q^2-1)(q-1)}-q^4\frac{(q^{n-2}-1)(q^{n-3}-1)}{(q^2-1)(q-1)}}
\ln\begin{bmatrix} n \\ 2 \end{bmatrix}_q\]
\[=4\frac{{(q^n-1)(q^{n-1}-1)}}
{{(q^n-1)(q^{n-1}-1)}-q^4(q^{n-2}-1)(q^{n-3}-1)}
\ln\begin{bmatrix} n \\ 2 \end{bmatrix}_q\]
\[=4\frac{{q^{2n-1}-q^n-q^{n-1}+1}}
{q^{n+2}+q^{n+1}-q^n-q^{n-1}-q^4+1}
\ln\begin{bmatrix} n \\ 2 \end{bmatrix}_q\sim4q^{n-3}(2n-4)\ln q.\]

The $G_q(n,2)$ graph is a strongly regular graph with parameters
\[\left(\frac{(q^n-1)(q^{n-1}-1)}{(q^2-1)(q-1)} , \frac{q(q+1)(q^{n-2}-1)}{q-1} ,
\frac{q^{n-1}-1}{q-1}+q^2-2 , (q+1)^2\right),\]
thus we can use Theorem \ref{UpperBoundAlmostStronglyRegular}, where
$N=\max\{a,b-1\}=a=\frac{q^{n-1}-1}{q-1}+q^2-2$, thus
\[\mu(G_q(n,2))\le\left\lceil \frac{v\ln\frac v{\sqrt2}}{k-N}\right\rceil\sim
\frac{q^{2n-4}\ln\frac{q^{2n-4}}{\sqrt2}}{q^{n-1}-2q^{n-2}}\sim
q^{n-3}(2n-4)\ln q.\]
We give two upper bounds obtained by construction.

\begin{definition}
	A line spread in $\PG(n-1,q)$ is a set of lines such that each point is covered by exactly
	one line of the set.
\end{definition}

\begin{theorem}[\cite{KissSzonyiBook} Theorem 5.14.] There is a line spread in $\PG(n-1,q)$
	if and only if $2\mid n$.
\end{theorem}

\begin{theorem}
	If $n$ is even, then the union of a line spread and a line pencil in $\PG(n-1,q)$ forms
	a resolving set of the graph $G_q(n,2)$, and therefore
\[\mu(G_q(n,2))\le\frac{q^{n-1}-1}{q-1}+\frac{q^{n}-1}{q^2-1}-1.\]
\end{theorem}

\begin{proof}
	Let $[P]$ be a pencil with carrier $P$, $\mathcal S^1$ be a spread and $S=[P]\cup\mathcal S^1$.
	Let $l_1$ and $l_2$ be two different lines in $\PG(n-1,q)$.

	\textbf{Case 1:} $l_1\cap l_2=\emptyset$. If $P\in l_1$, that is, $l_1\in[P]$ then
	$P\not\in l_2$, and thus $d(l_1,l_1)\ne d(l_1,l_2)$.
	If $P\not\in l_1$ and $P\not\in l_2$, then let us consider the plane spanned by $P$ and
	$l_1$. The line $l_2$ intersects this plane in at most one point, thus there are at least $q$
	lines through $P$ in this plane that intersect $l_1$ and do not intersect $l_2$.

	\textbf{Case 2:} $l_1\cap l_2\ne\emptyset$. If $l_1$ or $l_2$ is an element of $\mathcal S^1$
	then $d(l_1,l_1)\ne d(l_1,l_2)$. If neither $l_1$ nor $l_2$ is in $\mathcal S^1$, then
	let $P_1$ and $P_2$ be two points on $l_1\setminus l_2$ and let $a_1$ and $a_2$ be the lines of
	$\mathcal S^1$ through $P_1$ and $P_2$. Since $a_1$ and $a_2$ do not intersect each other,
	they cannot both lie in the plane spanned by $l_1$ and $l_2$, thus at least one of them
	does not intersect $l_2$.
\end{proof}

\begin{notation}
We will denote the set of lines through the point $P$ incident to the subspace
$\alpha$ by $[P|\alpha]$. That is,
\[[P|\alpha]:=\{l\in\mathcal L : P\in l , l\subseteq\alpha\},\]
\end{notation}

\begin{theorem}
	Let $n\ge5$ and let $P_1, P_2, \dots, P_n$ be points in general position in $\PG(n-1,q)$,
	and for each $i=1,\dots,n$ let $\alpha_i$ be the hyperplane spanned by the points
	$\{P_1,\dots,P_n\}\setminus\{P_i\}$. Then the set
	\[S:=[P_2|\alpha_1]\cup[P_3|\alpha_1]\ \cup\ [P_3|\alpha_2]\cup[P_4|\alpha_2]\ \cup\dots\cup\]
	\[[P_n|\alpha_{n-1}]\cup[P_1|\alpha_{n-1}]\ \cup\ [P_1|\alpha_n]\cup[P_2|\alpha_n]\]
	is a resolving set of $G_q(n,2)$, and therefore
	\[\mu(G_q(n,2))\le2n\frac{q^{n-2}-1}{q-1}.\]
\end{theorem}

\begin{proof}
	Let $a$ and $b$ be two different lines. We can assume that $a,b\notin S$.
	Here subscripts of the $P_i$ and $\alpha_i$ are taken modulo $n$.
		
	\textbf{Case 1:} There exists an $i\in\{1, \dots, n\}$ such that $a\subset\alpha_i$ and
	$b\not\subset\alpha_i$. At least one of $P_{i+1}$ and $P_{i+2}$ differs from
	$b\cap\alpha_i$. In its pencil inside $\alpha_i$ there are at least $q$
	lines that intersect $a$ but not $b$.

	\textbf{Case 2:} For each $i\in\{1, \dots, n\}$ both $a$ and $b$ either lie on $\alpha_i$
	or intersect $\alpha_i$. The subspace spanned by $a$ and $b$ is at least two-dimensional.

	\textbf{Case 2.1:} $a\cap b\ne\emptyset$. Let $Q=a\cap b$. 

	\textbf{Case 2.1.1:} There is exactly one hyperplane in $\{\alpha_1,\dots,\alpha_n\}$ that
	does not contain $Q$. Since the intersection of the other $n-1$ ones is one point of $\{P_1,\dots,P_n\}$,
	we can assume that $Q=P_1$. Then $a$ and $b$ intersect $\alpha_1$. Let $A_1=a\cap\alpha_1$ and
	$B_1=b\cap\alpha_1$.
	
	If $P_2A_1\ne P_2B_1$ or $P_3A_1\ne P_3B_1$, then $[P_2|\alpha_1]$ or $[P_3|\alpha_1]$
	resolves $\{a,b\}$. Otherwise $A_1, B_1, P_2$ and $P_3$ are collinear. Thus $a$ and $b$ are
	in the plane spanned by $P_1, P_2$ and $P_3$. It follows that $a$ and $b$ lie on $\alpha_{n-1}$.
	Then $\{a,b\}\subset[P_1|\alpha_{n-1}]\subset S$.

	\textbf{Case 2.1.2:} There are at least two hyperplanes in $\{\alpha_1,\dots,\alpha_n\}$ that
	do not contain $Q$, $\alpha_i$ and $\alpha_j$. If none of $[P_{i+1}|\alpha_i]$, $[P_{i+2}|\alpha_i]$
	$[P_{j+1}|\alpha_j]$ or $[P_{j+2}|\alpha_j]$ resolves $\{a,b\}$ then $P_{i+1},P_{i+2}, A_i$, and $B_i$ are
	collinear, and so are $P_{j+1},P_{j+2}, A_j$, and $B_j$. By the Veblen-Young axiom the lines
	$P_{i+1}P_{i+2}$ and $P_{j+1}P_{j+2}$ intersect. We can assume that the intersection point is
	$P_{i+2}=P_{j+1}$. Then the plane spanned by $a$ and $b$ is equal to the plane spanned by $P_{i+1}, P_{i+2}$,
	and $P_{i+3}$, and this plane is in $\alpha_i$, thus $Q\in\alpha_i$, a contradiction.

	\textbf{Case 2.2:} $a\cap b=\emptyset$. We can assume that $a$ and $b$ intersect $\alpha_1$.
	If neither $[P_2|\alpha_1]$ nor $[P_3|\alpha_1]$ resolves $\{a,b\}$ then $P_2,P_3, A_1$, and
	$B_1$ are collinear. The lines $a$ and $b$ cannot both lie on the plane spanned by $P_1, P_2$
	and $P_3$. Thus $[P_1|\alpha_n]$ resolves $\{a,b\}$.
\end{proof}

\begin{theorem} Assume that $n\ge4$. If $q\ge29$, then
	\[\mu(G_q(n,2))\ge\frac{3(q-1)q^{n-3}}{(q+1)}.\]
 If $q\le28$, then
	\[\mu(G_q(n,2))\ge\frac{\sqrt{2q} q^{n-3}}{e}.\]
\end{theorem}

\begin{proof}
We will denote the number of vertices and the degree of a vertex by $v$ and $d$. Then
\[v=\frac{q^n-1}{q-1}\ \frac{q^{n-1}-1}{q-1}\ \frac1{q+1}=
\frac{\sum_{i=0}^{n-1}q^i\sum_{i=0}^{n-2}q^i}{q+1}\]
\[\ge\frac{(q^{n-1}+q^{n-2})(q^{n-2}+q^{n-3})}{q+1}=\frac{q^{n-2}(q+1)q^{n-3}(q+1)}{q+1}\ge
q^{2n-4}.\]
The degree of a vertex is
\[d=\frac{q(q+1)(q^{n-2}-1)}{q-1}.\]
We give a lower bound for the quotient of $v$ and $d$:
\[\frac vd=\frac{(q^n-1)(q^{n-1}-1)}{q(q+1)(q^{n-2}-1)(q^2-1)}\ge
\frac{(q^n-1)(q^{n-1}-q)}{q(q+1)(q^{n-2}-1)(q^2-1)}\ge\]
\[\frac{q^n-1}{(q+1)^2(q-1)}=\frac{\sum_{i=0}^{n-1}q^i}{(q+1)^2}\ge
\frac{q^{n-1}+q^{n-2}}{(q+1)^2}=\frac{q^{n-2}(q+1)}{(q+1)^2}=\frac{q^{n-2}}{(q+1)}.\]

By Lemma \ref{MainPropositionGeneralGraphs}
\[\mu(G_q(n,2))\ge
\min\left\{\frac{sv(u+1)}{d}\ ,
\ \frac{u\sqrt[u]{(1-s)v}}{e\sqrt[u]{2}}\ ,
\ \frac{(1-s)v}2
\right\}\]
Let $s=\frac{q-1}q$, $u=2$ and let us use the lower bounds for $v$ and $\frac vd$.
\[\mu(G_q(n,2))\ge\min\left\{3\frac{q-1}q \frac{q^{n-2}}{(q+1)}\ ,
\ \frac{2\sqrt{\frac1q q^{2n-4}}}{e\sqrt{2}}\ ,
\ \frac{q^{2n-4}}{2q}
\right\}=\]
\[=\min\left\{\frac{3(q-1)q^{n-3}}{(q+1)}\ ,
\ \frac{\sqrt{2q} q^{n-3}}{e}\ ,
\ \frac{q^{2n-5}}{2}
\right\}.
\]
It is easy to see that $\frac{\sqrt{2q} q^{n-3}}{e}<\frac{q^{2n-5}}{2}$ and thus
\[\mu(G_q(n,2))\ge q^{n-3}\min\left\{\frac{3(q-1)}{(q+1)}\ ,
\ \frac{\sqrt{2q}}{e}\right\}.\]
Since $q\ge2$ is an integer, the inequality $\frac{3(q-1)}{(q+1)}\le \frac{\sqrt{2q}}e$ holds if 
and only if $q\ge29$.
\end{proof}

\section{Moore graphs}
\begin{definition}
	A Moore graph of diameter $d$ is a $k$-regular graph whose girth (the length of its
	shortest cycle) is more than $2d$.
\end{definition}

In this section we consider the Moore graphs with diameter 2. These are $k$-regular graphs
with $k^2+1$ vertices. In \cite{MooreGraphsDiam2Paper} Hoffman and Singleton showed that
the possible degrees are 2, 3, 7 and 57. The 7-regular one is called the \emph{Hoffman-Singleton graph}.
The existence of the 57-regular one is unsolved.

Bonato et al. proved the following bounds.
\begin{theorem}[\cite{LocNumMetDimGrDiam2} Theorem 3.1.]
	Let $G_k$ be a $k$-regular Moore graph of diameter 2 and girth 5. If $k\ge3$ then
	\[k\le\mu(G_k)\le2k-3\]
\end{theorem}
Thus $\mu(G_7)\ge7$ and $\mu(G_{57})\ge57$.
Applying Theorem \ref{SimpleLowerBoundDiamTwo} we have the following:
\begin{theorem}
	Let $G_k$ be a $k$-regular Moore graph of diameter 2 and girth 5. Then
	\[\mu(G_k)\ge\frac{2k^2}{k+4}.\]
\end{theorem}
Thus $\mu(G_7)\ge9$, and $\mu(G_{57})\ge107$.
However in \cite{MetDimSmallDisRegStrReg} Bailey proved by computer calculation that
$\mu(G_7)=11$.

\section{Bilinear forms spaces}

\begin{definition}[\cite{BrouwerCohenNeumaierDistRegGraphBook}, Section 9.5.A]
	The \emph{bilinear forms graph} $H_q(n,d)$ is an induced subgraph of the Grassmann
	graph $G_q(n+d,d)$ whose vertices are the $d$-dimensional subspaces intersecting a
	fixed $n$-dimensional subspace trivially.
\end{definition}

The diameter of $H_q(n,d)$ is $\min\{n,d\}$. In this section we consider the $H_q(n,2)$
graphs. $H_q(n,2)$ is a strongly regular graph with parameters
\[\left(q^{2n}\ ,\ (q+1)(q^n-1)\ ,\ q^n+q^2-q-2\ ,\ q(q+1)\right).\]

Feng and Wang proved the following upper bound:
\begin{theorem}[\cite{FengWangBilFormGrMetric} Theorem 1.1.]
	If $n\ge d\ge2$, then
	\[\mu(H_q(n,d))\le q^{n+d-1+\left\lfloor\frac{d+1}n\right\rfloor}.\]
\end{theorem}
As a corollary,
\[\mu(H_q(2,2))\le q^{4}\ ,\ \ \ \ \mu(H_q(3,2))\le q^{5},\]
\[\textrm{if } n\ge4 \textrm{, then } \mu(H_q(n,2))\le q^{n+1}.\]

To get a new upper bound, let us apply Theorem \ref{UpperBoundAlmostStronglyRegular}. 
Since $a=q^n+q^2-q-2\ge b-1=q(q+1)-1$ if $n>2$ or $q>2$, we get the following theorem.
\begin{theorem} If $n>2$ or $q>2$, then
	\[\mu(H_q(n,2))\le\left\lceil\frac{q^{2n}\ln\frac{q^{2n}}{\sqrt2}}{q^{n+1}-q^2+1}\right\rceil\sim2nq^{n-1}\ln q.\]
\end{theorem}

Now we use Lemma \ref{MainPropositionGeneralGraphs} to get a lower bound:
\begin{theorem}
If $q\ge29$, then
	\[\mu(H_q(n,2))\ge3\frac{(q-1)q^{n-1}}{q+1},\]
and if $q\le28$, then
	\[\mu(H_q(n,2))\ge\frac{\sqrt{2q}q^{n-1}}{e}.\]
\end{theorem}

\begin{proof}
In Lemma \ref{MainPropositionGeneralGraphs}, let $u=2$ and $s=\frac{q-1}q$. Then
\[
	\mu(H_q(n,2))\ge
		\min\left\{\frac{sv(u+1)}{k}\ ,
		\ \frac{u\sqrt[u]{(1-s)v}}{e\sqrt[u]{2}}\ ,
		\ \frac{(1-s)v}2
		\right\}\]
	\[=\min\left\{\frac{3(q-1)q^{2n}}{q(q^n-1)(q+1)}\ ,
	\ \frac{2\sqrt{\frac {q^{2n}}q}}{e\sqrt{2}}\ ,
	\ \frac{q^{2n}}{2q}
	\right\}\ge
	q^{n-1}\min\left\{\frac{3(q-1)}{(q+1)}\ ,
	\ \frac{\sqrt{2q}}{e}\ ,
	\ \frac{q^n}{2}
	\right\}.\]	
	It is easy to see that $\frac{\sqrt{2q}}{e}\le\frac{q^n}{2}$, and
	$\frac{3(q-1)}{(q+1)}\le\frac{\sqrt{2q}}{e}$, if and only if $q\ge29$.
\end{proof}

\end{document}